\input ./style/arxiv-general.cfg
\documentclass[aos,MSNbibl,seceqn,dvips]{arximspdf}
\makeatletter
   \@ifpackageloaded{graphicx}{}{\usepackage{graphicx}}
\makeatother

\doi{10.1214/15-AOS1270REJ}
\volume{43}
\issue{4}
\pubyear{2015}
\firstpage{1463}
\lastpage{1470}
\docsubty{FLA}
\referstodoi{10.1214/14-AOS1270}

\makeatletter
\newcommand{\eqref}[1]{(\ref{#1})}
\newtheorem{theorem}{Theorem}[section]
\makeatother

\begin{document}
\begin{frontmatter}
\vspace*{12pt}
\title{Rejoinder to discussions of ``Frequentist coverage of adaptive
nonparametric
Bayesian credible sets''}
\runtitle{Rejoinder}

\begin{aug}
\author[A]{\fnms{Botond}~\snm{Szab\'{o}}\ead[label=e1]{b.szabo@tue.nl}},
\author[B]{\fnms{A. W.}~\snm{van der Vaart}\corref{}\ead[label=e2]{a.w.van.der.vaart@umail.leidenuniv.nl}}
\and
\author[C]{\fnms{J. H.}~\snm{van Zanten}\ead[label=e3]{hvzanten@uva.nl}}
\runauthor{B. Szab\'{o}, A. W. van der Vaart and J.  H. van Zanten}
\affiliation{TU Eindhoven, Leiden University and University of Amsterdam}
\address[A]{B. Szab\'o\\
Department of Mathematics\\
Eindhoven University of Technology\\
P.O. Box 513\\
5600 MB Eindhoven\\
The Netherlands\\
\printead{e1}}
\address[B]{A. W. van der Vaart\\
Mathematical Institute\\
Leiden University\\
P.O. Box 9512\\
2300 RA Leiden\\
The Netherlands\\
\printead{e2}}
\address[C]{J. H. van Zanten\\
Korteweg-de Vries Institute for Mathematics\\
University of Amsterdam\\
P.O. Box 94248\\
1090 GE Amsterdam\\
The Netherlands\\
\printead{e3}}
\end{aug}

%
\received{\smonth{1} \syear{2015}}


\end{frontmatter}

We thank the discussants for their supportive comments and interesting
observations. Many questions are still open and not all methodological
or philosophical
questions may have an answer. Our reply addresses only a subset of questions
and is organized by topic. A final section reviews recent work.


\section{Hierarchical Bayes credible sets}

Our paper considers empirical Bayes tuning of the posterior
distribution, whereas many Bayesians
might prefer to use a hierarchical Bayes approach. Ghosal and Rousseau
ask whether, or conjecture
that, the hierarchical Bayes procedure behaves similarly as the
likelihood based empirical
Bayes procedure. Indeed, we can show exactly the same coverage of
hierarchical Bayes credible sets for polished tail
truths. A counterexample showing that hierarchical Bayes credible sets
also do \emph{not} cover without
some restriction was already given in \cite{SzVZ3}, while the size of
such sets
follows from \cite{KSzVZ}. Thus, within the context of our paper there
is no difference between
the two schemes.

In the hierarchical Bayes approach we endow the regularity
hyperparameter $\alpha$ with a hyperprior
distribution $\lambda$, and then apply an ordinary Bayes method with
the overall prior, for some
upper bound $A$ (possibly dependent on $n$),
\[
\Pi(\cdot)=\int_0^A \Pi_\alpha(\cdot)
\lambda(\alpha)\,d\alpha.
\]
For $\Pi (\cdot | X^{(n)})$ the posterior distribution relative
to this prior,
a hierarchical Bayes credible ball centered at the posterior mean $\hat
\theta_n$ is defined by its
radius $\hat r_{n,\gamma}$:
%
\begin{equation}
\label{eqradius}
\Pi \bigl(\theta\dvtx  \|\theta- \hat\theta_{n}
\|_2 \le\hat{r}_{n,
\gamma}  |  X^{(n)} \bigr) = 1-
\gamma.
\end{equation}
We blow this up a bit, and for $L > 0$ consider
%
\begin{equation}
\label{EqCredibleSet}
\hat C_n(L) = \bigl\{\theta\in\ell^2\dvtx  \|
\theta- \hat\theta_{n}\| _2 \le L\hat {r}_{n, \gamma}
\bigr\}.
\end{equation}
Under a mild regularity condition on $\lambda$, similar to that in
\cite
{KSzVZ}, these sets cover polish tail truths.

\begin{theorem}
Suppose that there exist $c_1,c_2 \geq0$, $c_3$ and $c_4,c_5>0$,
with $c_3>1$ if $c_2=0$, such that
$c_4^{-1}\alpha^{-c_3} \exp(-c_2\alpha) \leq\lambda(\alpha) \leq
c_4\alpha^{-c_3} \exp
(-c_2\alpha)$, for all $\alpha>c_1$ and $\lambda(\alpha)\geq c_5$
for all $0<\alpha\leq
c_1$. Then for any positive $A,L_0, N_0$ there exists a constant
$L$ such that
%
\begin{equation}\label{eqPostCov}
\inf_{\theta_0\in\Theta_{pt}(L_0)} P_{\theta_0} \bigl(\theta_0\in
\hat C_n(L) \bigr)\rightarrow1.
\end{equation}
Furthermore,\vspace*{1pt} for $A=A_n\leq\sqrt{\log n}/(4\sqrt{\log\rho
\vee e})$ this is true with
a slowly varying sequence [$L:=L_n\lesssim(3\rho^{3(1+2p)})^{A
_n}$ works].
\end{theorem}

\begin{pf}
The probability of interest $\mathrm{P}_{\theta_0} (\|\hat\theta
_{n}-\theta
_0\|_2\leq L \hat{r}_{n,\gamma} )$ is bounded below by
\begin{eqnarray*}
&&\mathrm{P}_{\theta_0} \bigl(\|\theta_0-{\mathrm{E}}_{\theta_0}
\hat \theta_{n,\underline\alpha
_n}\|_2+\|\hat\theta_{n,\underline\alpha_n}-{
\mathrm{E}}_{\theta
_0}\hat\theta _{n,\underline\alpha_n}\|_2+\|\hat
\theta_{n}-\hat\theta _{n,\underline\alpha_n}\| _2\leq L
\hat{r}_{n,\gamma} \bigr).
\end{eqnarray*}
Therefore, the theorem follows from Theorem~5.1 of \cite{SzVZ2}, if
%
\begin{eqnarray}
\|\theta_0-{\mathrm{E}}_{\theta_0}\hat\theta_{n,\underline\alpha
_n}
\|_2&\lesssim &  n^{-\underline\alpha_n/(1+2\underline\alpha_n+2p)},
\nonumber
\\
P_{\theta_0} \bigl(\|\hat\theta_{n,\underline\alpha_n}-{\mathrm {E}}_{\theta_0}
\hat \theta_{n,\underline\alpha_n}\|_2\leq   C_1 n^{-\underline\alpha_n/(1+2\underline\alpha_n+2p)}
\bigr)&\rightarrow & 1,
\nonumber
\\
P_{\theta_0} \bigl(\|\hat\theta_{n}-\hat\theta_{n,\underline\alpha
_n}\|
_2\leq C_2 n^{-\underline\alpha_n/(1+2\underline\alpha_n+2p)} \bigr) &\rightarrow & 1,\label{eqDiffPostMean}
\\
P_{\theta_0} \bigl(\hat{r}_{n,\gamma}\geq  C_3
n^{-\overline\alpha
_n/(1+2\overline\alpha_n+2p)} \bigr)&\rightarrow & 1.\label{eqHBradius}
\end{eqnarray}
The first two assertions follow immediately from $(5.8)$ and $(5.9)$ of
\cite{SzVZ2}.

For the proof of $\eqref{eqHBradius}$, we first note that, for any
given $C_3>0$,
\begin{eqnarray*}
&& \Pi\bigl(\theta\dvtx  \| \theta-\hat\theta_{n}\|_2<C_3
n^{-{\overline\alpha _n}/({1+2\overline\alpha_n+2p})} | X^{(n)}\bigr)
\\
&&\qquad = \int_{\underline\alpha_n}^{\overline\alpha_n}\Pi_{\alpha
}\bigl(\theta\dvtx
\|\theta-\hat \theta_{n}\|_2 <C_3
n^{-{\overline\alpha_n}/({1+2\overline\alpha_n+2p})} |  X^{(n)}\bigr)  \lambda\bigl(\alpha |
 X^{(n)}\bigr)\,d\alpha\\
 &&\qquad\quad{}+o_{P_{\theta_0}}(1).
\end{eqnarray*}
The right side becomes bigger if we replace $\hat\theta_n$ by $\hat
\theta_{n,\alpha}$, as the latter is the
center of the Gaussian distribution $\Pi_\alpha(\cdot | X^{(n)})$, and
again bigger if we replace
$\overline\alpha_n$ in the rate inside the probability by $\alpha$.
From the proof of $(5.7)$ of \cite{SzVZ2}, it follows that there exists
a constant $C_3$, such that for every $\alpha$,
\begin{eqnarray*}
&&\Pi_\alpha \bigl(\theta\dvtx  \|\theta-\hat\theta_{n,\alpha}
\|_2\leq C_3 n^{-{\alpha}/({1+2\alpha+2p})} | X^{(n)} \bigr)
\leq1-2\gamma.
\end{eqnarray*}
Then the integral in the preceding display is asymptotically smaller
than $1-2\gamma$, whence
\eqref{eqHBradius} follows by the definition of $\hat r_{n,\gamma}$.

To prove $\eqref{eqDiffPostMean}$, we proceed similarly to the proof
of $(4.4)$ of \cite{SzVZ3}.
By Jensen's inequality,
%
\begin{eqnarray}
 \|\hat\theta_{n}- \hat\theta_{n,\underline\alpha_n} \|_2^2
&\leq & \int\|\hat\theta_{n,\alpha}- \hat
\theta_{n,\underline\alpha
_n}\|_2^2  \lambda\bigl(
\alpha|X^{(n)}\bigr)\,d\alpha
\nonumber
\\
\label{eqhelp2}
&\leq & \sup_{\alpha\in[\underline\alpha_n,\overline\alpha_n]}\| \hat\theta _{n,\underline\alpha_n}- \hat
\theta_{n,\alpha}\|_2^2 \\
&&{}+ \sup_{\alpha\notin[\underline\alpha_n,\overline\alpha_n]}
\| \hat\theta_{n,\underline\alpha_n}-\hat\theta_{n,\alpha}\| _2^2
\int_{\alpha\notin
[\underline\alpha_n,\overline\alpha_n]} \lambda\bigl(\alpha | X^{(n)}\bigr)\,d
\alpha.
\nonumber
\end{eqnarray}
We separately bound the two terms on the right side. First, as
$\underline\alpha_n\in[\underline\alpha_n,\overline\alpha_n]$,
by several applications of the triangle inequality,
\begin{eqnarray*}
&&\sup_{\alpha\in[\underline\alpha_n,\overline\alpha_n]}\| \hat \theta_{n,\alpha}- \hat
\theta_{n,\underline\alpha_n}\|_2 
\\
&& \qquad \leq 2 \sup_{\alpha\in[\underline\alpha_n,\overline\alpha_n]}\| {\mathrm{E}}_{\theta_0} \hat
\theta_{n,\alpha}- \theta_0\|_2 +2\sup
_{\alpha\in[\underline\alpha_n,\overline\alpha_n]}\| \hat \theta_{n,\alpha}- {\mathrm{E}}
_{\theta_0} \hat\theta_{n,\alpha}\|_2.
\end{eqnarray*}
As a consequence\vspace*{1pt} of $(5.8)$ and $(5.9)$ of \cite{SzVZ2}, this is
bounded above
by a multiple of $n^{-\underline\alpha_n/(1+2\underline\alpha
_n+2p)}$, with
$P_{\theta_0}$-probability tending to one.
For the second term, we first note that similar to the preceding display,
with $P_{\theta_0}$-probability tending to one,
\begin{eqnarray*}
&& \sup_{\alpha}\| \hat\theta_{n,\alpha}- \hat
\theta_{n,\underline
\alpha_n}\|_2\leq 2\sup_{\alpha}\|{
\mathrm{E}}_{\theta_0} \hat\theta_{n,\alpha}- \theta_0
\|_2 +2\sup_{\alpha}\| \hat\theta_{n,\alpha}- {
\mathrm{E}}_{\theta_0} \hat\theta_{n,\alpha}\|_2.
\end{eqnarray*}
As a consequence of $(5.10)$ and $(5.11)$ of \cite{SzVZ2}, this is
uniformly bounded
by a constant times $\|\theta_0\|_2^2+1\lesssim1$,
with $P_{\theta_0}$-probability tending to one. Furthermore, in view of
Section~7 of \cite{KSzVZ},
\begin{eqnarray*}
&&{\mathrm{E}}_{\theta_0} \int_{\alpha\notin[\underline\alpha
_n,\overline\alpha_n]} \lambda\bigl(\alpha
 | X^{(n)}\bigr)\,d\alpha \\
 &&\qquad\leq 2e^{-({C_4n^{{1}/({1+2\overline\alpha
_n+2p)}}})/({1+2\overline\alpha
_n+2p})}(\log
n)^{C_5}e^{C_6e^{\sqrt{\log n}/3}}
\\
&&\qquad\leq  2e^{-({C_4n^{{1}/({1+2\overline\alpha
_n+2p})}})/({2(1+2\overline\alpha_n+2p)})}\\
&&\qquad \lesssim   n^{-({2\overline\alpha_n})/({1+2\overline\alpha_n+2p})}.
\end{eqnarray*}
Therefore, by Markov's inequality, the second term on the right-hand side
of $\eqref{eqhelp2}$ is bounded above by a multiple of
$n^{-({2\overline\alpha_n})/({1+2\overline\alpha_n+2p})}$, which
is smaller
than the same rate at $\underline\alpha_n$,
with $P_{\theta_0}$-probability tending to one.
\end{pf}

For the adaptive size we note that similar to the proof of assertion
(4.5) of \cite{SzVZ3},
it can be shown that there exists a positive constant $C_7$ such that
\begin{eqnarray*}
&& P_{\theta_0}\bigl(\hat{r}_{n,\gamma}\leq C_7
n^{-\underline\alpha
_n/(1+2\underline\alpha_n+2p)}\bigr)\rightarrow1.
\end{eqnarray*}
Then following \cite{SzVZ2}, we get the rate adaptive size for Sobolev
balls, hyperrectangles, analytic balls,$\,\ldots$ etc.

\section{Shape of the credible sets: Bands versus balls}

All discussants pointed out that $L_2$-confidence sets are harder to
visualize than
confidence bands, that is, $L_\infty$-balls. We fully agree. See our
remarks on plotting below.

We chose to consider $L_2$-balls because
they fit naturally in our inverse problem setup and can be studied
theoretically with reasonable
ease. At the same time, we believe that they provide an accurate (or at
least not misleading)
rendering of the general phenomena surrounding adaptive credible sets.
We fully agree that it is of interest to work out
similar results for other norms and other situations.\vadjust{\goodbreak}

One of the theoretical difficulties to handle credible bands is to
describe the $L_\infty$-norm in terms of quantities
that are controllable under the prior and posterior. Several authors
(starting with \cite{GineNickl})
have recently obtained contraction rates in
this norm, and their work may well be extendible to adaptive credible sets.

Ghosal proves a rate of contraction for the uniform norm, for
parameters such that
$\sum_ii^\alpha|\theta_i|<\infty$ 
and a prior that depends on $\alpha$. He next argues heuristically that
the resulting credible sets
with an \emph{adaptive} choice of $\alpha$ will cover relative to the
uniform norm. This is possible,
but the particular empirical Bayes $\hat\alpha$ from our paper may for
many true parameters
not estimate Ghosal's $\alpha$, but a different value.

We have encountered similar phenomena when deriving contraction rates
and credible intervals for
(not necessarily continuous) linear functionals of the parameter; see
\cite{SzLinear}. Since point evaluations are linear
functionals, such credible intervals can be glued together into
$L_{\infty}$-credible bands, where
due to the Gaussianity one would expect at most a logarithmic factor to
be necessary to pass from
pointwise to simultaneous intervals. A difficulty is that Sobolev
regularity is not the most useful concept when
estimating a function at a point; one would like to employ a H\"older
norm. As a worst
case one loses a $1/2$ when passing from Sobolev to H\"older, and this
loss was seen to be real
for the minimax contraction rate in \cite{Bartek,KSzVZ}. The
likelihood-based empirical Bayes method seems to ``estimate'' the
Sobolev regularity of the truth.
In \cite{SzLinear} we have shown that coverage can be retained by
subtracting $1/2$ from the
estimate, thus under-smoothing the empirical Bayes posterior
distribution. In forthcoming
work with Sniekers, we note that the ordinary empirical Bayes procedure
may still give good
coverage for many true parameters, the loss of $1/2$ being really a
worst case comparison of the
two norms and coverage being connected to more subtle properties of the
true parameter.

\section{Simulation and plotting}
We included some pictures in the paper, and we feel that they nicely
illustrate the
limitations and strengths of adaptive Bayesian credible sets. The
pictures consist of individual plots of all
curves in the 95\% out of 2000 curves simulated from the posterior
distribution that are closest to the posterior
mean. Within the resolution of the pictures these curves form a ragged
grey band and it is tempting
to view this as a confidence band.

We may have mislead the reader to think that the pictures show the
$L_2$-credible set
that we study theoretically in the paper. However, as already noted,
$L_2$-balls are difficult to plot. To relate our plots to
these balls, it seems one would have to ``visually compute'' the
$L_2$-distance of the plotted curves
to the center of the band (the posterior mean), take the maximum distance,
and compare this to the $L_2$-distance of a tentative function to this
center, in order to see
whether this function is in the ball. This is hard to do. The pictures
are not formal credible bands either.
Still, they manage to give an impression of where the posterior
distribution puts its mass.

Low and Ma describe this difficulty very accurately. In particular, our
choice of making
exactly 2000 draws was rather arbitrary and, indeed, at other places we
have also produced
pictures showing just 20 draws (without rejecting any). All these
pictures seem to \emph{illustrate}
the effect of the bias--variance trade-off, and its possible failure, on
credible sets reasonably well.

Low and Ma also suggest a method for constructing $L_{\infty
}$-confidence bands from the $L_2$-credible
balls and apply it to the adaptive posterior distribution.
Bayesians will be delighted to see that the empirical Bayes method
performs satisfactorily in
their simulation study. The new concept of coverage introduced by these
authors, together
with Cai, is interesting.

Castillo also addresses the discrepancy between our analytic definition
of a credible ball and
our small simulation study. He points out that the radius can be
simulated more precisely.
He also suggests that simulating curves from distributions that
are rougher than the posterior might be useful to fill out the gap
between the support
of the posterior and the ball. This is an interesting suggestion, but
we would be reluctant to simulate
from other distributions than the posterior distribution. We imagine
that this
could be queried in many ways, for example, to produce bands,
intervals for specific
functionals or perhaps even of qualitative aspects of parameters, but
we would support the
Bayesian view that the posterior distribution gives a full report of
the analysis.

Nickl and Castillo \cite{CasNickl2} have introduced an approach toward
credible sets based on a nonparametric
Bernstein--von Mises theorem. Nickl writes to be ``unsure to which
extent $\ell_2$-credible balls
are applied in current practice as claimed in the introduction of (our
paper),'' and next suggests
that ``Practitioners may prefer ($\ldots$) to compute credible balls in
$\mathbb{H}$-spaces.'' Castillo
wonders about our opinion that ``no method that avoids dealing with the
bias--variance trade-off will
properly quantify the uncertainty.'' We do not believe we have claimed
that $\ell_2$-\emph{balls}
are routine in practice; if we did, then we retract that claim here. We
do claim that posterior
distributions are routinely used for uncertainty quantification, often
by simulating from it.
Then a main finding of 50 years of
nonparametric statistics, theory and practice, is that the
bias--variance trade-off drives
everything, setting it apart from classical, parametric statistics,
which deals mostly with variance, as bias
is negligible, particularly in the large-sample limit. The
$\ell_2$-setting of our paper incorporates the bias--variance trade-off,
and hence we believe that our theoretical
results are relevant. It appears that Nickl and Castillo's
``Bernstein--von Mises theorem in $\mathbb{H}$-space'' removes
bias, essentially by parameterizing the function as a collection of
smooth functionals that
can be estimated as the parameters in classical parametric models, with
neglect of bias. Their work
is very intriguing and pretty. However, as it explains away bias, we
found it difficult to believe that it solves
the nonparametric problem. It is still more intriguing that pictures by
Ray in \cite{Ray2}, which are based
on the $\mathbb{H}$-spaces of \cite{CasNickl2}, look similar to ours.
Possibly that is because
these pictures do not show their suggested set, just as our pictures
are deficient in this sense.
This deserves further investigation.

\section{Other priors}

The discussants pose the question whether our results extend to other priors
than the $N(0, i^{-1-2\alpha})$-priors in our paper. We believe the answer
is affirmative: it appears
that the polished tail condition is not linked to the form of the priors.

One reason to believe this are preliminary results, of ourselves and in
a forthcoming thesis of
Sniekers at Leiden University, about priors of the form
\[
\Pi_{\tau}=\prod_{i=1}^{\infty}N
\bigl(0,\tau^2i^{-1-2\alpha}\bigr),
\]
where $\alpha$ is fixed, but $\tau$ is adapted to the data, by either an
empirical or hierarchical Bayes
method. For empirical Bayes we plug the marginal maximum likelihood
estimator $\hat\tau_n$ of $\tau$ into
the posterior distributions for fixed $\tau$, and construct adaptive
credible sets of the form
%
\begin{equation}
\label{EqCredibleSet} \hat C_n(L) = \bigl\{\theta\in\ell^2\dvtx  \bigl\|
\theta- \hat\theta _{n}(\hat\tau_n)\bigr\| _2 \le
Lr_{n, \gamma}(\hat\tau_n) \bigr\},
\end{equation}
where $\hat\theta_n(\tau)$ is the posterior mean and $r_{n,\gamma
}(\tau)$ satisfies
%
\begin{equation}
\label{eqradius} \Pi \bigl(\theta\dvtx  \bigl\|\theta- \hat\theta_{n}(\tau)
\bigr\|_2 \le r_{n,
\gamma}(\tau)  | X^{(n)} \bigr) = 1-
\gamma.
\end{equation}

\begin{theorem}\label{ThmScaling}
For any $A, L_0, N_0$ there exists a constant $L$ such that
\begin{eqnarray*}
&& \inf_{\theta_0\in\Theta_{pt}(L_0)} P_{\theta_0} \bigl(\theta_0\in
\hat C_n(L) \bigr)\rightarrow1.
\end{eqnarray*}
\end{theorem}

Hierarchical Bayes credible sets will similarly cover. However,
these sets have the disadvantage
that they may be unnecessarily big. In our paper \cite{SzVZ} we proved
that the corresponding posterior distributions
contract at the minimax rate over Sobolev balls of regularity $\beta
<\alpha
+1/2$, but only at the
suboptimal rate $n^{-(1+2\alpha)/(4+4\alpha)}$ if $\beta>\alpha
+1/2$. The latter
suboptimal rate is partially due to the variance of the posterior
distribution, and hence,
in the case that $\beta>\alpha+1/2$, the radius of the credible balls will
be suboptimal as well.

\section{Choice of basis}

Rousseau and Castillo point out that the polished tail condition is
dependent on the chosen basis,
whereas one might hope or expect the set of ``good behaving'' true parameters
not to depend on the basis.

In inverse problems the eigenbasis of the operator $K^*K$
plays, implicitly or explicitly, an important role to
describe the problem \cite{Bartek,Bartek2,Ray,Agapiou} and, hence, it
is natural to assume the polished tail
condition with respect to this basis. Other bases were explored
in recent work \cite{KnapikSalomond}, but
a good link between the operator and the prior seems always needed.

In ``direct problems'' one can consider any basis. This then determines
both the
prior and the polished tail condition. The prior, or rather collection
of priors, will
be chosen to model a scale of models that is thought to capture the
true parameter.
In our situation these were Sobolev spaces, which are naturally
described in a basis.
That the polished tail condition will adopt the same basis seems not unnatural.
After all, ``good-behaving'' is not an absolute property of a
parameter, but is relative
to a method, which is the one induced by the prior in this case.

There is a good scope for extensions to other models and priors.
In our case the coefficients could be modeled differently than
independent and Gaussian,
although both seem natural. We imagine that similar results as in our
paper can easily be written
down for double-indexed bases, as wavelets, thus moving closer to the earliest
works on self-similarity. More challenging will be priors such as
Dirichlet mixtures, which
are known to adapt to the bandwidth in the (normal) kernel. What can be said
about their coverage?

\section{Further references}
The paper \cite{Ray2} derives an adaptive and nonparametric version of
the Bernstein--von Mises
theorem, using techniques developed in \cite{CasNickl2} and \cite{KSzVZ},
under a self-similarity restriction, and next applies this result to
construct adaptive credible
sets. The same paper also considers spike and slab priors and
$L_{\infty
}$-credible bands.
The author of \cite{Belitser} investigates credible sets from an oracle
perspective. He considers truncated (finite dimensional)\vadjust{\goodbreak}
Gaussian priors and shows that the empirical Bayes approach chooses the
optimal truncation level
under a (slightly) extended version of the polished tail condition.
This family of priors is
relatively wide and contains a member that attains the minimax
posterior contraction rate
for every regularity class $S^{\beta}$. The authors of \cite
{SniekersvdV} have followed up their work
with investigating adaptive pointwise credible sets using rescaled
(integrated) Brownian motion
as a prior in the nonparametric regression model.
Random smoothing spline priors with Gaussian weights on the spline
coefficients are shown in \cite{Serra} to give
honest credible sets in the nonparametric regression problem under the
self-similarity condition.

\printaddresses
\end{document}